\input amstex.tex
\documentstyle{amsppt}
\overfullrule=0pt
\pageheight{19cm}

\leftheadtext{C.D.Hill and M.Nacinovich}

\magnification=\magstep1

\font\sc=cmcsc10

\hyphenation{pseu-do-con-ca-ve}
\hyphenation{s-o-n-d-e-r-f-o-r-s-c-h-u-n-g-s-b-e-r-e-i-c-h}
\hyphenation{de-pen-den-ce}
\newcount\q
\newcount\x
\newcount\t

\long\def\se#1{\advance\q by 1
\x=0  \t=0 \bigskip
\noindent
\S\number\q \quad
{\bf {#1}}\par
\nopagebreak}

\long\def\thm#1{\advance\x by 1
\bigskip\noindent%
{\sc Theorem \number\q.\number\x}
\quad{\sl #1}
\par 
\noindent}

\long\def\prop#1{\advance\x by 1
\bigskip\noindent%
{\sc Proposition \number\q.\number\x}
\quad{\sl #1} 
\noindent}

\long\def\thml#1#2{\advance\x by 1
\bigskip\noindent
{\sc Theorem \number\q.\number\x ({#1})}
\quad{\sl #2} 
\noindent}

\long\def\lem#1{\advance\x by 1
\medskip\noindent
{\sc Lemma \number\q.\number\x}
\quad{\sl #1} 
\noindent}

\long\def\form#1{\global\advance\t by 1
$${#1} \tag \number\q.\number\t$$}
\long\def\cor#1{\advance\x by 1
\bigskip\noindent%
{\sc Corollary \number\q.\number\x}%
\quad{\sl #1} 
\noindent}

\def\dimo{\demo{Proof}
}

\hyphenation{dep-end-en-ce}
\topmatter
\title
Fields of $CR$ meromorphic functions
\endtitle
\author
C.Denson Hill and Mauro Nacinovich
\endauthor
\address C.Denson Hill - Department of Mathematics, SUNY at Stony Brook,
Stony Brook NY 11794, USA \endaddress
\email dhill\@math.sunysb.edu \endemail
\address Mauro Nacinovich - Dipartimento di Matematica -
Universit\`a di Roma Tor Vergata -
via della Ricerca Scientifica, 00133 - Roma,  Italy \endaddress
\email nacinovi\@mat.uniroma2.it\endemail
\toc
\widestnumber\subhead{3.1}
\specialhead {} Introduction \endspecialhead
\head 1. Preliminaries\endhead
\head 2. $CR$ meromorphic functions on compact $CR$ manifolds\endhead
\head 3. Analytic and algebraic dependence of $CR$ meromorphic functions
\endhead
\head 4. The field of $CR$ meromorphic functions\endhead
\head 5. The Chow theorem for $CR$ manifolds\endhead
\head 6. Projective embedding \endhead
\endtoc
\abstract
Let $M$ be a smooth compact $CR$ manifold of $CR$ dimension
$n$ and $CR$ codimension $k$, which has a certain local extension
property $E$. In particular, if $M$ is pseudoconcave, it has property
$E$. Then the field $\Cal K(M)$ of $CR$ 
meromorphic functions on $M$ has
transcendence degree $d$, with $d\leq n+k$. If
$f_1,\, f_2,\, \hdots ,\, f_d$ 
is a maximal set of algebraically independent
$CR$ meromorphic functions on $M$, then $\Cal K(M)$ is a simple finite
algebraic extension of the field 
$\Bbb C(f_1,\,f_2,\,\hdots,\,f_d)$ of rational
functions of the $f_1,\,f_2,\,\hdots ,\,f_d$. 
When $M$ has a projective embedding,
there is an analogue of Chow's theorem, and $\Cal K(M)$ is isomorphic
to the field $\Cal R(Y)$ of rational functions on an irreducible 
projective algebraic variety $Y$, and $M$ has a $CR$ embedding in
$\roman{reg}\,Y$. The equivalence between algebraic dependence and analytic
dependence fails when condition $E$ is dropped.
\endabstract
\endtopmatter

\document
\centerline{{\sc Introduction}}\smallskip
In a beautiful paper Siegel [Si], improving upon an idea of Serre [Se],
managed to give simple proofs of the basic theorems concerning
algebraic dependence and transcendence degree for the field of meromorphic
functions on an arbitrary compact complex manifold; thereby generalizing
classical results about the field of Abelian functions on a complex $n$
dimensional torus. For a detailed discussion of the now nearly 150 year
history of these matters, see the paper of Siegel. His proofs were based
on his extension to $n$ dimensions of the classical Schwarz lemma.
Later, following almost exactly Siegel's argument, Andreotti and
Grauert [AG] were able to show that the Siegel modular group, which plays
a pivotal role in the study of algebraic fields of automorphic functions,
is pseudoconcave. Later Andreotti [A] generalized these kinds of results
to general pseudoconcave complex manifolds and spaces; again following
Siegel's method.\par
In the present work we replace the compact complex manifold of Siegel
by a smooth compact pseudoconcave $CR$ manifold $M$ of general
$CR$ dimension $n$ and $CR$ codimension $k$, and study algebraic
dependence, transcendence degree and related matters for the field
$\Cal K(M)$ of $CR$ meromorphic functions on $M$. Again we follow
the method of Siegel, based on the Schwarz lemma, and we incorporate
some ideas used by Andreotti. Actually we are able to obtain 
results under a condition on the $CR$ manifold that is weaker
than pseudoconcavity, which we call condition $E$. In particular
we obtain an analogue of Chow's theorem [C] for compact $CR$ 
manifolds. In the situation where $M$ has a projective embedding,
we are able to identify $\Cal K(M)$ with the field
$\frak R(Y)$ of rational functions on an irreducible algebraic
variety $Y$, in which $M$ has a generic $CR$ embedding that avoids
the singularities of $Y$. We show that the possibility for $M$
to have a projective embedding is equivalent to the existence of 
a complex $CR$ line bundle over $M$ having certain properties. In this
context, it is interesting to note that the general abstract notion of
a complex $CR$ line bundle $F$ over a $CR$ manifold is such that
$F$ may fail to be locally $CR$ trivializable, even in the case where
$M$ is $CR$ embeddable [HN8]. \par
For more information about pseudoconcave $CR$ manifolds, we refer the
reader to the foundational paper [HN3], to the many examples in
[HN8], and to [HN1], [HN2], $\hdots$, [HN11], as well as [BHN], [DCN],
and [L].

\se{Preliminaries}
An abstract smooth almost $CR$ manifold
of type $(n,k)$ consists of: a connected smooth paracompact manifold
$M$ of dimension $2n+k$, a smooth subbundle $HM$ of $TM$ of rank
$2n$, that we call the {\it holomorphic tangent space} of $M$,
and a smooth complex structure $J$ on the fibers of $HM$.\par Let
$T^{0,1}M$ be the complex subbundle of the complexification $\Bbb CHM$
of $HM$, which corresponds to the $-\sqrt{-1}$ eigenspace of $J$:
\form{T^{0,1}M\,=\,\{X+\sqrt{-1}JX\, \big{|}\, X\in HM\}\, .} We say
that $M$ is a $CR$ manifold if, moreover, the formal integrability
condition \form{\left[ \Cal C^\infty(M,T^{0,1}M),\,
\Cal C^\infty(M,T^{0,1}M)\right]\,
\subset\, \Cal C^\infty(M,T^{0,1}M)\, } 
\noindent
holds. 
When $k=0$, via the Newlander-Nirenberg theorem, we recover the
definition of a complex manifold.
\par
Next we define ${T^*}^{1,0}M$ as the
annihilator of $T^{0,1}M$ in the complexified cotangent bundle $\Bbb C
T^*M$.  We denote by $Q^{0,1}M$ the quotient bundle $\Bbb C T^* M /
{T^*}^{1,0}M$, with projection $\pi_Q$.  It is a rank $n$ complex vector
bundle on $M$, dual to $T^{0,1}M$.  The $\bar\partial_M$--operator
acting on smooth functions is defined by $\bar\partial_M= \pi_Q\circ d$.
A local trivialization of the bundle $Q^{0,1}M$ on an open set $U$ in
$M$ defines $n$ smooth sections $\bar L_1$, $\bar L_2$, $\hdots$, $\bar
L_n$ of $T^{0,1}M$ in $U$; hence
\form{\bar\partial_M u\, = \,
\left(\bar L_1u, \bar L_2u,\hdots, \bar L_nu\right)\, ,}
where $u$ is a

function in $U$.  Solutions $u$ of $\bar\partial_Mu=0$ are called $CR$
functions. We denote by $\Cal C\Cal R(U)$ the space of {\it smooth 
($\Cal C^\infty$) functions}
on an open subset $U$ of $M$ that satisfy $\bar\partial_Mu=0$.
Note that $\bar\partial_M$ is a homogeneous first order partial
differential operator and hence the space $\Cal C\Cal R(U)$
is a commutative algebra with respect to the multiplication of
functions.
We denote by $\Cal C\Cal R_{M,a}=\varinjlim_{U\ni a}
\Cal C\Cal R(U)$ the local ring of 
germs of smooth $\Cal C\Cal R$ functions at
$a\in M$.
\par
\medskip
Let $M_1$, $M_2$ be two abstract 
smooth $CR$ manifolds, with holomorphic tangent spaces
$HM_1$, $HM_2$, and partial complex structures $J_1$, $J_2$, respectively.
A smooth map $f:M_1@>>>M_2$ is $CR$ if $f_*(HM_1)\subset HM_2$,
and $f_*(J_1v)=J_2f_*(v)$ for every $v\in HM_1$.\par
A $CR$ {\it embedding}
\footnote{In this case
we shall often identify $M$ with the submanifold $\phi(M)$ of $X$.}
$\phi$ of an abstract $CR$ manifold $M$ into
a complex manifold $X$, with complex structure
$J_X$, is a $CR$ map which is a smooth embedding
satisfying $\phi_*(H_aM)=\phi_*(T_aM)\cap J_X(\phi_*(T_aM))$
for every $a\in M$. We say
that the embedding is {\it generic} if the complex dimension of $X$
is $(n+k)$, where $(n,k)$ is the type of $M$. 
\par
Let $M$ be a smooth abstract $CR$ manifold of type $(n,k)$.
We say that $M$ is {\it locally embeddable} at $a\in M$, if $a$ 
has an open neighborhood $\omega_a$ in $M$ which admits a $CR$ embedding
into some complex manifold $X_a$. In this case we can always take
for $X_a$ an open subset of $\Bbb C^{n+k}$ and assume that the embedding
$\omega_a\hookrightarrow X_a$ is generic. The property of being
locally embeddable at $a$ is equivalent to the fact that
there exist an open neighborhood $\omega_a$ of $a$ and functions
$f_1,f_2,\hdots,f_{n+k}\in\Cal C\Cal R(\omega_a)$ such that
\form{d\,f_1(a)\wedge d\,f_2(a)\wedge\cdots \wedge d\,f_{n+k}(a)
\neq 0\, .}
The functions $f_1,f_2,\hdots,f_{n+k}$ can be taken to be
the restrictions to $\omega_a$ of the coordinate functions
$z_1,z_2,\hdots,z_{n+k}$ of $X_a\subset\Bbb C^{n+k}$. For this
reason one can say that they provide {\it $CR$ coordinates} on $M$
near $a$.
\bigskip
The {\it characteristic bundle} $H^0M$ is defined to be the
annihilator of $HM$ in $T^*M$. Its purpose it to parametrize
the Levi form: recall that the {\it Levi form}
of $M$ at $x$ is defined for $\xi\in H^0_xM$ and $X\in H_xM$ by
\form{\Cal L(\xi;X)\, = \, d\tilde\xi(X,JX)=\langle\xi,
[J\tilde X,\tilde X]\rangle\, ,}
where $\tilde\xi\in\Cal C^\infty(M,H^0M)$ and $\tilde
X\in\Cal C^\infty(M,HM)$
are smooth extensions of $\xi$ and $X$.
For each fixed $\xi$ it is a Hermitian quadratic form for the
complex structure $J_x$ on $H_xM$.\par
A $CR$ manifold $M$ is said to be {\it $q$-pseudoconcave} if 
the Levi
form $\Cal L(\xi;\,\cdot\,)$ has at least $q$ negative and $q$
positive eigenvalues for every 
$a\in M$ and every nonzero $\xi\in H_a^0M$.
\par
By the term {\it pseudoconcave $CR$ manifold} $M$ we mean an abstract
$CR$ manifold which is: ($i$) locally embeddable at each point, and
($ii$) $1$-pseudoconcave.
\smallskip
In this paper we shall be concerned with $CR$ manifolds
$M$ of type $(n,k)$ which have a certain property $E$ ($E$ is for
{\it extension}). $M$ {\sl is said to have property} $E$ iff there is
an $E$-pair $(M,X)$. By an $E$-pair we mean that
\roster
\item"($i$)"\quad $M$ is a generic $CR$ submanifold of the complex
manifold $X$, and
\item"($ii$)"\quad for each $a\in M$, the restriction map induces
an isomorphism $\Cal O_{X,a} @>>>\Cal C\Cal R_{M,a}$.
\endroster
\smallskip
\noindent
{\sc Remark}\quad {\sl If $M$ is a pseudoconcave $CR$ manifold,
then $M$ has property $E$.}\par
In fact, property ($i$) for a pseudoconcave $M$ was proved in
Proposition 3.1 of [HN3]; however, Theorem 1.3 below gives a new
simplified proof of this fact. Property ($ii$) for a pseudoconcave
$M$ was proved in [BP], [NV]; however, a very short proof of this
fact is also given by Theorem 13.2 in [HN7]. Thus property $E$ is to
be regarded as a somewhat weaker hypothesis on $M$ than pseudoconcavity.
\smallskip
When $k=0$, so $M$ is of type $(n,0)$, then $M$ is an $n$-dimensional  
complex manifold, and we obtain an $E$ pair by choosing $X=M$. Hence
we adopt the convention that {\it any complex manifold has property}
$E$.\par
When $n=0$, so $M$ is of type $(0,k)$, then $M$ is a smooth totally
real $k$-dimensional manifold, and we can never obtain an $E$-pair,
(unless $M=X=\text{a point}$), because then any smooth function belongs to
$\Cal C\Cal R(M)$.
\thm{Let $(M,X)$ be an $E$-pair. Then for any open set $\omega\subset M$
there is a corresponding open set $\Omega\subset X$ such that
\roster
\item"($i$)"\quad $\Omega\cap M=\omega$, and
\item"($ii$)"\quad $r:\Cal O(\Omega)@>>>\Cal C\Cal R(\omega)$ is
an isomorphism.
\endroster} \edef\teoaa{\number\q.\number\x}
\dimo
We fix a Hermitian metric $g$ on $X$, with associated distance $d(x,y)$.
Let $a\in\omega$ and consider
\form{\Cal F_n=\{(f,\tilde f)\in\Cal C\Cal R(\omega)\times 
\Cal O(B(a,\frac{1}{n}))\, | \, \tilde f=f\;\text{on}\;B(a,\frac{1}{n})
\cap\omega\,\}\, .} 
Here $B(a,\frac{1}{n})$ denotes the ball of radius $\frac{1}{n}$ in $X$,
centered at $a$. Note that each $\Cal F_n$ is a closed subspace of
a Fr\'echet-Schwartz space, and hence a Fr\'echet-Schwartz space
itself. For each $n$, the map
\form{\pi_n:\Cal F_n\ni (f,\tilde f) @>>> f\in \Cal C\Cal R(\omega)}
is linear and continuous. By our hypothesis,
\form{\bigcup_{n=1}^\infty{\pi_n(\Cal F_n)}=\Cal C\Cal R(\omega)\, .}
Hence by the Baire category theorem, some $\pi_{n_0}(\Cal F_{n_0})$
is of the second category. It follows from a theorem of Banach that
$\pi_{n_0}:\Cal F_{n_0}@>>>\Cal C\Cal R(\omega)$ is surjective.
Now we denote $B(a,\frac{1}{n})$ by $B_a$.\par
Next we fix a tubular neighborhood $U$ of $M$ in $X$, with
$\pi:U@>>>M$ denoting the orthogonal projection. By letting
$\rho\in\Cal C^\infty(\omega,\Bbb R^+)$ vary, we produce a fundamental
system of open neighborhoods
\form{\Omega_\rho=\{z\in U\, | \, \pi(z)\in\omega\, ,\;\text{and}\;
d(z,\pi(z))<\rho(\pi(z))\,\}}
of $\omega$ in $X$. We choose $\rho_0\in\Cal C^\infty(\omega,\Bbb R^+) $ 
such that
\form{\Omega_{\rho_0}\subset\bigcup_{a\in\omega}{B_a}\, .}
Since $M$ is a deformation retract of $\Omega_{\rho_0}$, and the local
holomorphic extension of $CR$ functions from the generic $M$ is unique,
the different extensions to each $B_a$ of a given $f\in\Cal C\Cal R(\omega)$
match at points of $\Omega_{\rho_0}$. This completes the proof with
$\Omega=\Omega_{\rho_0}$.
\enddemo
\smallskip
\cor{In the situation of Theorem \teoaa\; we have, in addition, that
\roster
\item"($iii$)" \quad 
If $f\in\Cal C\Cal R(\omega)$,
and $f$ vanishes of infinite order at $a\in\omega$, then $f\equiv 0$
in the connected component of $a$ in $\omega$.
\item"($iv$)" \quad $(r^*f)(\Omega)=f(\omega)$.
\item"($v$)"\quad If $|f|$ has a local maximum at a point $a\in\omega$, then
$f$ is constant on the connected component of $a$ in $\omega$.
\endroster}
\dimo
If $f\in\Cal C\Cal R(\omega)$ vanishes of infinite order at $a\in\omega$,
then also
$r^*f$ vanishes of infinite order at $a$ and, by the strong
unique continuation of holomorphic functions, $r^*f$ vanishes
identically in the connected component of $\Omega$
containing $a$, and we obtain ($iii$).\par
To prove ($iv$), we assume by contradiction that $r^*f$ takes some value
$z_0\in\Bbb C$ at some point of $\Omega$, but that $f$ does not assume
that value at any point of $\omega$. Then the function
\form{g\, = \, \dsize\frac{1}{f-z_0}}
belongs to $\Cal C\Cal R(\omega)$, and has no holomorphic extension to
$\Omega$, contradicting ($ii$).\par
By 
($ii$) and ($iv$), a local maximum of $f\in\Cal C\Cal R(\omega)$ at
$a\in\omega$, is also a local maximum of $r^*f$ at $a\in\Omega$;
thereby $r^*f$ is constant on the connected component of $a$ in $\Omega$
and we obtain ($v$).
\enddemo
\cor{Let $(M,X)$ and $(N,Y)$ be $E$-pairs, and let $f:M@>>>N$ be a
smooth $CR$ isomorphism. Then there are $E$-pairs $(M,X')$ and
$(N,Y')$, with $X'\subset X$ and $Y'\subset Y$, such that $f$
extends to a biholomorphic diffeomorphism 
$\tilde f:X'@>>>Y'$.} \edef\teoac{\number\q.\number\x}
\dimo
We first consider the case where $X$ and $Y$ are open sets in $\Bbb C^{n+k}$.
By Theorem \teoaa, there is an open neighborhood $\Omega$ of $M$ in $X$
where $f$ has a holomorphic extension $\tilde f:\Omega@>>>\Bbb C^{n+k}$.
By shrinking $\Omega$ to $\Omega'$, we can arrange that 
$\tilde f(\Omega')\subset Y$ and the Jacobian determinant of $\tilde f$
is different from zero on $\Omega'$. Likewise there is an open neighborhood
$\Omega''$ of $N$ in $Y$ where $f^{-1}=g$ extends to $\tilde g$, with
$g(\Omega'')\subset\Omega'$ and the Jacobian determinant of $\tilde g$
being nonzero in $\Omega''$. By uniqueness of holomorphic extension of
$CR$ functions from $M$ to $\tilde g(\Omega'')$, it follows that
$\tilde g\circ\tilde f=\text{identity}$ on a neighborhood of $M$ in $X$.\par
Now we consider the general case. Introducing local holomorphic 
coordinates charts on $X$ and $Y$, we may use the special case above to
produce local holomorphic extensions. The local holomorphic extensions
patch together, by unique continuation, to give the desired $\tilde f$.
\enddemo
We may now use Corollary \teoac\; to show that $M$ having property $E$
is actually a local property of $M$.
\thm{$M$ has property $E$ if and only if for each $a\in M$, there is
an open neighborhood $\omega_a$ of $a$ in $M$ such that $\omega_a$
has property $E$.}
\dimo
By hypothesis we have an $E$-pair $(\omega_a,X_a)$, for each $a\in M$.
We can assume that $\omega_a\Subset M$, and that $\pi_a:X_a@>>>\omega_a$
is the orthogonal projection from a tubular neighborhood, with a
distance function $d_a(x,y)$. By Corollary \teoac, whenever
$\omega_a\cap\omega_b\neq\emptyset$, there are open neighborhoods
$X_{ab}$ of $\omega_a\cap\omega_b$ in $X_a$ and
$X_{ba}$ of $\omega_a\cap\omega_b$ in $X_b$,
and a unique biholomorphic map $\tilde f_{ab}:X_{ab}@>>>X_{ba}$,
extending the identity map on $\omega_a\cap\omega_b$.
We may select a locally finite open covering $\{\omega_a\}$ of $M$,
parametrized by $a\in A\subset M$. By shrinking, we refine the
$\{\omega_a\}$ to an open covering $\{\omega_a'\}$, with
$\omega'_a\Subset\omega_a$. With $\epsilon_a>0$ sufficiently small, we
define
\form{X'_a=\pi_a^{-1}(\omega_a')\cap\{d_a(x,\omega'_a)<\epsilon_a\}\, ,}
so as to have
\form{\pi_a^{-1}(\omega'_a\cap\omega'_b)\cap X'_a\subset X_{ab}\, ,}
for all $b\in A$ such that $\omega'_b\cap\omega'_a\neq\emptyset$.
Set
\form{X'_{ab}=\tilde f_{ab}^{-1}(X_{ab})\cap X'_a\, .}
Then $X$ is obtained by gluing together the $X'_a$'s, by
\form{\CD
X'_a\supset X'_{ab}@>{\tilde f_{ab}}>>X'_{ba}\subset X'_b\, .\endCD}
This completes the proof.
\enddemo
\medskip \def\CR{\Cal C\Cal R}
We now turn to the object of main concern in this paper, which are the
$CR$ meromorphic functions on an $M$ satisfying property $E$.
The ring $\CR(\omega)$ of smooth ($\Cal C^\infty$) $CR$ functions on
$\omega\subset M$ is an integral domain if $\omega$ is connected.
Let $\varDelta(\omega)$ be the subset of $\CR(\omega)$
of divisors of zero; i.e. $\varDelta(\omega)$ is the set
of those $CR$ functions on $\omega$ which vanish in some connected
component of $\omega$. Let $\Cal M(\omega)$ be the quotient ring of
$\CR(\omega)$ with respect to $\CR(\omega)\setminus\varDelta(\omega)$.
This means that $\Cal M(\omega)$ is the set of the equivalence classes
of pairs 
$(p,q)$
with $p\in\CR(\omega)$
and $q\in\CR(\omega)\setminus\varDelta(\omega)$.
The equivalence relation $(p,q)\sim (p',q')$ is defined by
$pq'=p'q$. If $\omega'\subset\omega$ is an inclusion of open subsets of $M$,
the restriction map $r^{\omega'}_{\omega}:\CR(\omega)@>>>\CR(\omega')$
sends $\CR(\omega)\setminus\varDelta(\omega)$
into $\CR(\omega')\setminus\varDelta(\omega')$ and thus induces a homomorphism
of rings:
\form{r^{\omega'}_{\omega}:\Cal M(\omega) @>>>\Cal M(\omega')\, .}
We obtain in this way a presheaf of rings. We shall call the corresponding
sheaf $\Cal M$ the {\it sheaf of $CR$ meromorphic functions on $M$}.
By a $CR$ meromorphic function on an open set $\omega\subset M$, we mean
a continuous section $f$ of $\Cal M$ over $\omega$. If $\omega$ is
connected, the space of all such sections $\Cal K(\omega)$ forms
a field. Since we always assume that $M$ is connected, we have in particular
$\Cal K(M)$, {\it the field of $CR$ meromorphic functions on $M$}.
\medskip
We recall these standard notions: Let $\Bbb F$ be a field and
$\Bbb F_0\subset\Bbb F$ a subfield. Then $f_1,f_2,\hdots,f_\ell\in\Bbb F$
are said to be {\it algebraically dependent} over $\Bbb F_0$ iff
there is a nonzero polynomial $P\in\Bbb F_0[x_1,x_2,\hdots,x_\ell]$ with
coefficients in $\Bbb F_0$ such that
\form{P(f_1,f_2,\hdots,f_\ell)=0\, ;}
otherwise they are called {\it algebraically independent}. The
{\it transcendence degree} of $\Bbb F$ over $\Bbb F_0$ is the
cardinality of a maximal set $S\subset\Bbb F$ such that every finite
subset of $S$ is algebraically independent over $\Bbb F_0$. If the
transcendence degree of $\Bbb F$ over $\Bbb F_0$ is zero, we say that
$\Bbb F$ is {\it algebraic} over(or is an {\it algebraic extension of})
$\Bbb F_0$.
The cardinal $[\Bbb F:\Bbb F_0]$ denotes the dimension of
$\Bbb F$ over $\Bbb F_0$, as a vector space. 
The field $\Bbb F$ is said to be
a {\it simple algebraic extension of} $\Bbb F_0$ if there exists
an element $\theta\in\Bbb F$ such that any $f\in\Bbb F$ can be 
written as a polynomial in $\theta$ with coefficients in $\Bbb F_0$.
When $\Bbb F_0$ has characteristic zero, the primitive element theorem
says that any finite algebraic extension of $\Bbb F_0$ is simple.
\medskip
Finally we discuss the general notion of a smooth {\it complex $CR$ line
bundle} $F@>\pi>>M$, which was introduced in [HN8]. By this we mean that
$F$ is a smooth complex line bundle over $M$ such that:
\roster
\item"($i$)"\quad $F$ and $M$ are smooth abstract $CR$ manifolds of type
$(n+1,k)$ and $(n,k)$, 
respectively,
\item"($ii$)"\quad $\pi:F@>>>M$ is a smooth $CR$ submersion,
\item"($iii$)"\quad 
$F\oplus F\ni(\xi_1,\xi_2)@>>>\xi_1+\xi_2\in F$ and
$\Bbb C\times F\ni(\lambda,\xi)@>>>\lambda\cdot\xi\in F$
are $CR$ maps.
\endroster
Note that the Whitney sums $F\oplus F$ and $\Bbb C\times F$ have natural
structures of smooth $CR$ manifolds of type $(n+2,k)$; see [HN8].
There we also introduced the notion of the tangential $CR$ operator
$\bar\partial^F_M$, acting on smooth sections of $F$. We may take
a smooth (not necessarily $CR$) trivialization 
$(U_\alpha,\sigma_\alpha)$ of $F$, where $\sigma_\alpha$ is a smooth 
non vanishing section of $F$ on $U_\alpha$. Then a smooth section $s$
of $F$ has a local representation $s=s_\alpha\sigma_\alpha$ in $U_\alpha$,
where $s_\alpha$ is a smooth complex valued function in $U_\alpha$,
and $s_\alpha=g_{\alpha\beta}s_\beta$ in $U_\alpha\cap U_\beta$, with
$g_{\alpha\beta}=s_\beta/s_\alpha$. In each $U_\alpha$ the tangential
$CR$ operator acting on $s$ has a representation of the form:
\form{\bar\partial^F_M s=(\bar\partial_Ms_\alpha+A_\alpha s_\alpha)\otimes
s_\alpha\, ,}\edef\formulaar{\number\q.\number\t}
where $A_\alpha\in\Cal C^\infty(U_\alpha,\Cal Q^{0,1}M)$ and
$\bar\partial_MA_\alpha=0$.
On $U_\alpha\cap U_\beta$ we have:
\form{A_\beta-A_\alpha=g_{\beta\alpha}\bar\partial_Mg_{\alpha\beta}\, ,
\;\text{with}\; g_{\beta\alpha}=g_{\alpha\beta}^{-1}\, .}
If $s$ satisfies $\bar\partial^F_Ms=0$, it is called a $CR$ section of
$F$.\par
The $\ell$-th tensor power $F^\ell$ of $F$ is still a smooth complex
$CR$ line bundle over $M$, which can be defined in the same trivialization,
and we have
\form{\bar\partial^{F^\ell}_Mt=(\bar\partial_Mt+\ell A_\alpha t)\otimes
\sigma_\alpha^\ell}\edef\formulaat{\number\q.\number\t}
in $U_\alpha$, where $t$ is a smooth section of $F^\ell$.
\par
If the local trivialization is a $CR$ trivialization, then the $(0,1)$ forms
$A_\alpha$ in (\formulaar) and (\formulaat) are equal to zero.
On the other hand, if the $\bar\partial_M$-closed forms $A_\alpha$ are
locally $\bar\partial_M$-exact, then $F$ and $F^\ell$ are locally $CR$
trivializable.
\smallskip
Let $f\in\Cal K(M)$, where we now assume that $M$ has property $E$. Then we
can associate to $f$ a smooth complex line bundle $F@>\pi>>M$, which is
locally $CR$ trivializable. By definition, $f$ has local representations
\form{f=\dsize\frac{p_a}{q_a}\quad\text{on}\quad \omega_a\, ,}
with $p_a,q_a\in\CR(\omega_a)$. Moreover we may arrange that their
holomorphic extensions $\tilde p_a$ and $\tilde q_a$ to a neighborhood
$\Omega_a$ of $\omega_a$ in $X$ have no nontrivial common factor
at each point of $\Omega_a$. Then there are uniquely determined
non vanishing functions $g_{ab}\in\CR(\omega_a\cap\omega_b)$ such that
\form{q_a=g_{ab}q_b\quad\text{on}\quad\omega_a\cap\omega_b\, .}
The $\{g_{ab}\}$ are then the transition functions of a smooth 
complex $CR$ line bundle $F$ over $M$, and $F$ is therefore locally
$CR$ trivializable. The $\{p_a\}$ and $\{q_a\}$ give global smooth
sections $p$ and $q$ of $F$ over $M$, whose quotient $p/q$ is the
$CR$ meromorphic function $f$.
\smallskip
Let us return now to the smooth complex $CR$ line bundle $F@>\pi>>M$,
which may not be locally $CR$ 
trivializable. In this context, it is natural
to consider smooth abstract $CR$ manifolds $M$, which may not have
property $E$, but which are {\it essentially pseudoconcave}, as defined
\footnote{$M$ is {\it essentially pseudoconcave} 
iff it is {\it minimal}, i.e. does not
contain germs of $CR$ manifolds with the same $CR$ dimension and
a smaller $CR$ codimension,
and admits a Hermitian
metric on $HM$ for which the traces of the Levi forms are zero at each
point.}
in [HN8]. The important consequence of the assumption of
essential pseudoconcavity on $M$ is that one has the weak unique
continuation property for $CR$ sections of $F$. Note that
$1$-pseudoconcave abstract $CR$ manifolds are essentially pseudoconcave.
Under these assumptions we can give a more general notion of what is
a $CR$ meromorphic function on $M$: We associate a $CR$ meromorphic
function $f$ to any pair $(p,q)$, where $p$ and $q$ are smooth global
$CR$ sections of a smooth complex $CR$ line bundle $F@>\pi>>M$, with
$q\not\equiv 0$. Another pair $(p',q')$, which are smooth
$CR$ global sections of another such $F'@>{\pi'}>>M$, with $q'\not\equiv 0$,
define the same $f$ iff $pq'=p'q$ as sections of $F\otimes F'$.
Note that $f=p/q$ is a well defined smooth $CR$ function where $q\neq 0$.
With this more general definition, we get a new collection
$\hat{\Cal K}(M)$ of objects called $CR$ meromorphic functions on $M$.
Observe that $\hat{\Cal K}(M)$ is a field. For an essentially
pseudoconcave $M$, which has property $E$, $\Cal K(M)$ is a subfield
of $\hat{\Cal K}(M)$. If in addition $M$ is $2$-pseudoconcave, then all
smooth complex $CR$ line bundles over $M$ are locally $CR$ trivializable,
and then $\Cal K(M)=\hat{\Cal K}(M)$.
\se{$CR$ meromorphic functions on compact $CR$ manifolds.}
Let $M$ be a connected smooth compact $CR$ manifold of type $(n,k)$,
having property $E$. Then:
\thm{The field $\Cal K(M)$ of $CR$ meromorphic functions on $M$ has
transcendence degree over $\Bbb C$ 
less or equal to $n+k$.}\edef\teoba{\number\q.\number\x}
\par
\quad Setting $k=0$ above, we recover Satz 1 in Siegel [Si].
\dimo
According to the discussion in \S 1, the statement means: Given
$n+k+1$ $CR$ meromorphic functions $f_0,f_1,\hdots,f_{n+k}$ on $M$,
there exists a non zero polynomial with complex coefficients
$F(x_0,x_1,\hdots,x_{n+k})$ such that
\form{F(f_0,f_1,\hdots,f_{n+k})\equiv 0\quad\text{on}
\quad M\, .}\edef\formba{\number\q.\number\t}
From the preceding section, we may regard $M$ as a generic $CR$
submanifold of an $n+k$ dimensional complex manifold $X$.\par
For each point $a\in M$ there is a connected open coordinate 
neighborhood $\Omega_a$, in which the holomorphic coordinate
$z_a$ is centered at $a$. We choose $\Omega_a$ in such a way that
$\omega_a=\Omega_a\cap M$ is a connected neighborhood of $a$ in $M$.
Moreover we can arrange that, for $j=0,1,\hdots,n+k$, each $f_j$
has a representation
\form{f_j=\dsize\frac{p_{ja}}{q_{ja}}\quad\text{on}\quad\omega_a}
with \edef\formubb{\number\q.\number\t}
$p_{ja}$ and $q_{ja}$ being smooth $CR$ functions in $\omega_a$.
According to Theorem \teoaa\; we may also assume that the restriction 
map $\Cal O(\Omega_a)@>>>\CR(\omega_a)$ is an isomorphism.
For each $CR$ function $g$ on $\omega_a$, we denote its unique
holomorphic extension 
to $\Omega_a$ by $\tilde{g}$.\edef\teoba{\number\q.\number\x}
By a careful choice of the $p_{ja}$ and $q_{ja}$, and an additional
shrinking of $\omega_a$, $\Omega_a$, we can also arrange that
\form{\tilde{f}_j=
\dsize\frac{\tilde{p}_{ja}}{\tilde{q}_{ja}}\quad\text{on}
\quad\Omega_a\, ,}\edef\formubc{\number\q.\number\t}
with the functions $\tilde{p}_{ja}$ and $\tilde{q}_{ja}$
being holomorphic and having no 
nontrivial common factor at each point in a neighborhood of
$\overline{\Omega}_a$. For each pair of points $a,b$ on $M$ we have the
transition functions
\form{\tilde{q}_{ja}=g_{jab}\tilde{q}_{jb}\, 
,}\edef\formubd{\number\q.\number\t}
which are holomorphic and non vanishing on a neighborhood of
$\overline{\Omega}_a\cap\overline{\Omega}_b$. Again, for each $a\in M$
we consider the polydiscs:
\form{K_a=\{|z_a|\leq r_a\}\quad\text{and}\quad L_a=\{|z_a|<e^{-1}r_a\}\, 
,}\edef\formbe{\number\q.\number\t}\noindent
where $|z_a|$ denotes the max norm in $\Bbb C^{n+k}$, and $r_a>0$ is
chosen so that $K_a\Subset\Omega_a$. By the compactness of $M$, we may
fix a finite number of points $a_1,a_2,\hdots,a_m$ on $M$, such that the
$L_{a_1},L_{a_2},\hdots,L_{a_m}$ provide an open covering of $M$. Then
we choose positive real numbers $\mu$ and $\nu$ to provide the bounds:
\form{|g_{0ab}|<e^\mu\qquad\text{and}\qquad \left|\dsize\prod_{j=1}^{n+k}{
g_{jab}}\right|<e^\nu} \edef\formbf{\number\q.\number\t}
on $\overline{\Omega}_a\cap\overline{\Omega}_b$ for
$a,b=a_1,a_2,\hdots,a_m$.\par
Consider a polynomial with complex coefficients to be determined later,
\hfill\linebreak
$F(x_0,x_1,\hdots,x_{n+k})$ of degree $s$ with respect to $x_0$ and
of degree $t$ with respect to each $x_i$ for $i=1,2,\hdots,n+k$. The number
of coefficients to be determined is
\form{A=(s+1)\cdot (t+1)^{n+k}\, .}
Now, letting $a$ stand for any one of the $a_1,a_2,\hdots,a_m$, we
introduce the functions
\form{Q_a=\tilde q^s_{0a}\prod_{j=1}^{n+k}{\tilde q_{ja}^t}\,,
\quad P_a=Q_aF(\tilde f_0,\tilde f_1,
\hdots,\tilde f_{n+k})}\edef\formbh{\number\q.\number\t}\noindent
which are holomorphic on a neighborhood of $\overline{\Omega}_a$.
For a positive integer $h$, to be made precise later, we wish to impose
the condition, for $a=a_1,a_2,\hdots,a_m$, that $P_a$ vanishes to order
$h$ at $a$. In terms of our local coordinates $z_a$, this means that
all partial derivatives of order $\leq h-1$ must vanish at $z_a=0$.
This imposes a certain number of linear homogeneous conditions on the
unknown coefficients of the polynomial $F$. The number of such conditions
is 
\form{B=m\binom{n+k+h-1}{n+k}\leq m\, h^{n+k}\, .}
If we can arrange that $B<A$, then this system of linear homogeneous
equations has a non trivial solution.
\par
However, in order to apply the Schwarz lemma later, we need also to
arrange that $s$, $t$ and $h$ satisfy
\form{\mu\,s\,+\,\nu\,t\, <\, h\, .}\edef\formbl{\number\q.\number\t}
To this end we fix $s$ to be an integer with $s>m\nu^{n+k}$. Thus, for 
each positive $h$, we denote by $t_h$ the largest positive integer 
satisfying
$st_h^{n+k}<mh^{n+k}$.
In this way we obtain that
\form{B\leq mh^{n+k}\leq s\left(t_h+1\right)^{n+k}<(s+1)
\left(t_h+1\right)^{n+k}=A\, .}
On the other hand, since $t_h@>>>\infty$ as $h@>>>\infty$, by choosing
$h$ sufficiently large we have
\form{m\left(\dsize\frac{\mu s}{t_h}+\nu\right)^{n+k}
<s\, ,}
which implies (\formbl) for $t=t_h$. Set
\form{\Upsilon=\max_{1\leq i\leq m}\,
\max_{K_{a_i}}{|P_{a_i}|}\, .}
This maximum is obtained at some point $z^*$ belonging to some
$K_{a^*}$, for $a^*$ equal to some one of $a_1,a_2,\hdots,a_m$.
Since $z^*\in K_{a^*}\subset\Omega_{a^*}$, because of our choices of
the $\omega_a$, $\Omega_a$, according to ($iv$) in Corollary 1.2,
there is another point $z^{**}\in \omega_{a^*}$ such that
\form{P_{a^*}(z^*)=P_{a^*}(z^{**})\, .}\edef\formbp{\number\q.\number\t}
But the point $z^{**}$ belongs to some $L_{a^{**}}\subset K_{a^{**}}$,
where $a^{**}$ is one of the $a_1$, $a_2$, $\hdots$,
$a_m$. Hence by the Schwartz
lemma of Siegel [Si] we obtain
\form{\left|P_{a^{**}}(z^{**})\right|
\leq \Upsilon\, e^{-h}\, .} \edef\formbq{\number\q.\number\t}
However
\form{P_{a^{*}}(z^{**})=P_{a^{**}}(z^{**})\left[
g_{0a^*a^{**}}^s(z^{**})\dsize\prod_{j=1}^{n+k}{g^t_{ja^*a^{**}}(z^{**})}
\right]\, .} 
Hence from (\formbf), (\formbp), (\formbq) we obtain
\form{\Upsilon=\left|P_{a^{*}}(z^{**})\right|\leq \Upsilon\, 
e^{\mu s+\nu t -h}\, .}
By (\formbl) this implies that $\Upsilon=0$. Hence each $P_{a_j}\equiv 0$,
which in turn yields $F(\tilde f_0,\tilde f_1,\hdots,\tilde f_{n+k})\equiv 0$.
Therefore restricting to $M$ we get (\formba). This completes the proof.
\enddemo
\se{Analytic and algebraic dependence of $CR$ meromorphic functions}
Let $f_0,f_1,\hdots,f_\ell\in \Cal K(M)$. We say that they are
{\it analytically dependent} if
\form{df_0\wedge df_1\wedge \cdots \wedge df_\ell=0\quad\text{
where it is defined.}}\edef\formca{\number\q.\number\t}
\thm{Let $M$ be a connected smooth compact $CR$ manifold of type
$(n,k)$, having property $E$. Let $f_0,f_1,\hdots,f_\ell\in \Cal K(M)$.
Then they are algebraically dependent over $\Bbb C$ if and only if they
are analytically dependent.}\edef\teorca{\number\q.\number\x}
\demo{Proof} 
First we observe that algebraic dependence implies analytic dependence.
Assume that there is a nontrivial polynomial $F$, with complex coefficients,
of minimal total degree, such that $F(f_0,f_1,\hdots,f_\ell)\equiv 0$.
Then
\form{\dsize\sum_{j=0}^\ell{
\dsize\frac{\partial F}{\partial x_j}(f_0,f_1,\hdots,f_\ell)\, 
df_j=0}}\edef\formcb{\number\q.\number\t}\noindent
where it is defined. It follows that some coefficient in (\formcb)
is a nonzero $CR$ meromorphic function on $M$. This implies (\formca)
on an open dense subset of $M$, and hence whenever it is defined.\par
For the proof in the other direction, we can assume that
$f_1,\hdots,f_\ell$ are analytically independent. Our task is to show
that there exists a nonzero polynomial with complex coefficients
$F(x_0,x_1,\hdots,x_\ell)$ such that
\form{F(f_0,f_1,\hdots,f_\ell)\equiv 0
\quad\text{on}\quad M\, .}\edef\formcc{\number\q.\number\t}\noindent
To this end we repeat the proof of Theorem \teoba, with $n+k$ replaced
by $\ell$, down to the line below (\formbh).
We shall replace $s,t,\nu,A,B$ by new $s',t',\nu',A',B'$.
After that we choose additional points $a_1',a_2',\hdots,a_m'$ with
$a'_j\in\omega_{a_j}$ and $a'_j$ sufficiently close to $a_j$, for 
$j=1,2,\hdots,m$. These points are chosen so that
\form{K_{a'_j}=\left\{\left|z_{a_j}-z_{a_j}(a_j')\right|\leq r_{a_j}\right\}
\Subset \Omega_{a_j}\, ,}\noindent
the 
$L_{a'_j}=\left\{\left|z_{a_j}-z_{a_j}(a_j')\right|< e^{-1}r_{a_j}\right\}$
still give an open covering of $M$, the
$\tilde f_0$, $\tilde f_1$,
$\hdots$, $\tilde f_\ell$ are holomorphic at each $a_j'$,
and $\tilde f_1,\hdots,\tilde f_\ell$ can be completed to a system
of holomorphic coordinates in a neighborhood of each $a'_j$. This is
possible because the set of points 
on $M$, where the $\tilde f_1,\hdots,\tilde f_\ell$
are holomorphic, and the $d\tilde f_1,\hdots, d\tilde f_\ell$ are linearly
independent,
is open and dense. Our assumption (\formca) that the
$f_0,f_1,\hdots,f_\ell$ are analytically dependent implies that, near each
point $a'_j$, $\tilde f_0$ is a holomorphic function of 
$\tilde f_1,\hdots,\tilde f_\ell$. We modify the proof of Theorem \teoba\;
by requiring that the holomorphic functions $P_{a_j}$ vanish to order
$h$ at $a'_j$, for $j=1,\hdots, m$. To accomplish this, we require that
$F(\tilde f_0,\tilde f_1,\hdots,\tilde f_\ell)$ vanish to order $h$
at each $a'_j$. This amount to requiring that all partial derivatives
of order $\leq (h-1)$ of $F(\tilde f_0,\tilde f_1,\hdots,\tilde f_\ell)$
with respect to $\tilde f_1,\hdots,\tilde f_\ell$ should vanish at $a'_j$.
The number of homogeneous linear equations is now
\form{B'=m\binom{\ell+h-1}{\ell}\leq m \, 
h^{\ell}\, .}\noindent
Now we fix an integer $s'$ with $s'>m(\nu')^\ell$.
Just as in the proof of Theorem \teoba, we choose $h$ sufficiently large,
and take $t'$ to be the largest positive integer satisfying
$s'(t')^\ell<mh^\ell$, so as to obtain
\form{B'\leq mh^\ell<(s'+1)(t'+1)^\ell
=A'\, ,}\edef\formcf{\number\q.\number\t}\noindent
\form{\mu s'+\nu' t'<h\, .}\noindent
By (\formcf) we can choose a nontrivial $F$, of degree $s'$
in $x_0$ and of degree $t'$ in each $x_1,\hdots,x_\ell$, such that
all $P_{a_j}$ vanish to order $h$ at $a'_j$.
Set
\form{\Upsilon=\max_{1\leq j\leq m}\, \max_{K_{a'_j}}{\left|P_{a_j}\right|\,
}.}\noindent
This maximum is attained at some point $z'$ belonging to some
$K_{a'_{j_0}}$. Since $z'\in K_{a'_{j_0}} \subset\Omega_{a_{j_0}}$,
as before by ($iv$) in Corollary 1.2 there is another point 
$z''\in\omega_{a_{j_0}}$ such that
\form{P_{a_{j_0}}(z')=P_{a_{j_0}}(z'')\, 
.}\noindent
This point $z''$ belongs to some $L_{a'_{j_1}}\subset K_{a'_{j_1}}$.
So by the Schwarz lemma we obtain
\form{\left| P_{a_{j_1}}(z'')\right|\leq \Upsilon\, e^{
\mu\,s'\,+\,\nu\,t'\,-\,h}\, .}
Thus as before we obtain
\form{\Upsilon=\left| P_{a_{j_1}}(z'')\right|\leq \Upsilon\, e^{
\mu\,s'\,+\,\nu\,t'\,-\,h}\, ,}
showing that $\Upsilon=0$. This implies (\formcc), completing the proof.
\enddemo
In order to make the exposition more clear, we have divided the discussion
into two parts; however Theorem \teoba\; is a direct consequence of
Theorem 3.1.
\smallskip
Algebraic dependence always implies analytic dependence. However, in
the absence of property $E$  , the converse may be false. We give a
general counterexample: 
\prop{Let $M$ be a connected smooth compact $CR$ manifold of type
$(n,k)$. Assume that $M$ has a smooth $CR$ immersion into some Stein
manifold. Then
\roster
\item condition $E$ is violated, and
\item there exists an infinite sequence of smooth $CR$ functions on $M$,
any two of which are analytically dependent, and which are algebraically
independent over $\Bbb C$.
\endroster}
\demo{Proof} By the embedding theorem for Stein manifolds, we can assume that
$M$ has a smooth $CR$ immersion in some $\Bbb C^N$. Then by a result
in [HN1], there exists a point $x_0$ in $M$ and a $\xi\in H^0_{x_0}M$
such that the Levi form $\Cal L_{x_0}(\xi,\,\cdot\,)$ is positive definite
on $H_{x_0}M$. This implies that a small neighborhood of $x_0$ in $M$ is
contained in the smooth boundary of a strictly pseudoconvex domain
$\Omega$ in $\Bbb C^N$. It is well known that there are holomorphic
functions in $\Omega$, which are smooth on $\overline\Omega$, and cannot
be holomorphically extended beyond $x_0$. Thus condition $E$ is violated
at $x_0$.\par
To prove (2) we first observe that some coordinate $z_1$ must be non constant
on $M$. Consider the sequence of holomorphic functions
$\left\{f_\lambda(z)=e^{z_1^\lambda}\left|\, \lambda\in\Bbb N\right.\right\}
\subset\Cal O(\Bbb C^N)$, and denote their pull-backs to $M$ by
$\{f^*_\lambda\}$.
Clearly any two of them are analytically dependent. Assume, by contradiction,
that some finite collection $f^*_{\lambda_1},
f^*_{\lambda_2},\hdots, f^*_{\lambda_m}$, are algebraically dependent.
Then there is a nontrivial polynomial $P$, with complex coefficients,
such that
\form{P(f^*_{\lambda_1},
f^*_{\lambda_2},\hdots, f^*_{\lambda_m})=0\quad\text{on}\quad M\, .}
Since $e^{z_1^{\lambda_1}},e^{z_1^{\lambda_2}},\hdots,
e^{z_1^{\lambda_m}}$ are algebraically independent in $\Cal O(\Bbb C)$,
the entire function 
$z_1 @>>> P\left(e^{z_1^{\lambda_1}},e^{z_1^{\lambda_2}},\hdots,
e^{z_1^{\lambda_m}}\right)$ is not constant and has isolated zeroes.
Because $M$ is connected, it follows  that the pullback on $M$ of
the function $z_1$ is constant, contradicting our assumption that
the coordinate $z_1$ was not constant on $M$.
\enddemo
\se{The field of $CR$ meromorphic functions}
Suppose $M$ is a connected compact $CR$ manifold of type $(n,k)$, having
property $E$. We have:
\thm{Let $d$ be the transcendence degree of  $\Cal K(M)$ over $\Bbb C$,
and let $f_1,f_2,\hdots, f_d$ be a maximal set of algebraically independent
$CR$ meromorphic functions in $\Cal K(M)$. Then $\Cal K(M)$ is a
simple finite algebraic extension of the field
$\Bbb C(f_1,f_2,\hdots, f_d)$ of rational functions of 
$f_1,f_2,\hdots, f_d$.}
Setting $k=0$ above, and taking the special case where $d=n$, we recover
Satz 2 of Siegel [Si].\par\edef\teorda{\number\q.\number\x}
The theorem is a consequence of the following:
\prop{Let $f_1,f_2,\hdots,f_\ell$ be $CR$ meromorphic functions in 
$\Cal K(M)$. Then there exists a positive integer 
$\kappa=\kappa(f_1,f_2,\hdots,f_\ell)$ such that every $f_0\in\Cal K(M)$,
which is algebraically dependent on $f_1,f_2,\hdots,f_\ell$, satisfies a
nontrivial polynomial equation of degree $\leq\kappa$ whose coefficients
are rational functions of $f_1,f_2,\hdots,f_\ell$.}
\demo{Proof}
Without any loss of generality we may assume that $f_1,f_2,\hdots,f_\ell$
are algebraically independent. By Theorem \teorca\; they are also
analytically independent. This puts us in the situation of the second
half of Theorem \teorca. The difference, however, is that we use
{\sl only} the functions $f_1,f_2,\hdots,f_\ell$ in (\formubb),
(\formubc) and (\formubd) to determine the $\omega_a$, $\Omega_a$.
In this way the numbers $m$ and $\nu'$ depend only on
$f_1,f_2,\hdots,f_\ell$. We fix the integer
$s'>m(\nu')^\ell$ as before. The proof of Theorem \teorca\; shows
that any $CR$ meromorphic function $f\in\Cal K(M)$, which can be 
represented on each $\omega_{a_i}$, for $i=1,2,\hdots,m$,
as the quotient $p_i/q_i$ of two $CR$ functions globally defined on
$\omega_{a_i}$, satisfies a nontrivial polynomial equation
of degree less or equal to $\kappa=s'$, with 
coefficients in $\Bbb C(f_1,f_2,\hdots,f_\ell)$. This reduces our
task to showing that $f_0$ has such a representation.\par
By hypothesis our given $f_0$ satisfies an equation
\form{G_0f_0^\lambda+G_1f_0^{\lambda-1}+\cdots+G_
\lambda=0}\edef\formuda{\number\q.\number\t}
\noindent
where $G_0,G_1,\hdots,G_\lambda$ are polynomials in 
$f_1,f_2,\hdots,f_\ell$, and
$G_0$ is not identically $0$ in $M$. Let $\sigma$ be an upper bound
for the degrees of the $G_0,G_1,\hdots,G_\lambda$, with respect to
each of the $f_1,f_2,\hdots,f_\ell$. We set
\form{\tilde Q_a=\dsize\prod_{j=1}^\ell{\tilde q_{ja}^\sigma}\, 
,}
\form{\tilde H_{a\alpha}=\tilde Q_a^\alpha\tilde G_0^{\alpha-1}
\tilde G_\alpha\qquad (\alpha=1,2,\hdots,\lambda)\, ,
}
\form{\tilde S_a=\tilde Q_a\tilde G_0\tilde f_0\, .
}
Multiplying (\formuda) by $\tilde Q_a^\lambda\tilde G_0^{\lambda-1}$,
we obtain that
\form{\tilde S^\lambda+\tilde H_{a1}\tilde S^{\lambda-1}+
\cdots +\tilde H_{a\lambda}=0\quad\text{on}\quad\Omega_a\, .
}\edef\formude{\number\q.\number\t}
Note that $\tilde Q_a\tilde G_\beta$ ($\beta=0,1,\hdots,\lambda$)
and the $\tilde H_{a\alpha}$ ($\alpha=1,2,\hdots,\lambda$) are
holomorphic functions on $\Omega_a$, and $\tilde S_a$ is meromorphic
on an open neighborhood of $\omega_a$ in $\Omega_a$. 
Since $\tilde S_a$ satisfies (\formude), it is locally bounded,
and hence actually holomorphic. Then the restrictions
\form{p_{0a}=\tilde S_a\left|_{\dsize\omega_a}\right.
\quad\text{and}\quad q_{0a}=\tilde Q_a
\tilde G_0\left|_{\dsize\omega_a}\right.
}
\noindent
are $CR$ functions on $\omega_a$, and
\form{f_0=\dsize\frac{p_{0a}}{q_{0a}}\quad\text{on}\quad \omega_a\, .
}
The proof is complete.
\enddemo\edef\teordb{\number\q.\number\x}
Now we explain what is the point of Theorem \teorda. Consider a
maximal set $f_1,f_2,\hdots, f_d$ of algebraically independent
$CR$ meromorphic functions on $M$, where $d$ is the transcendence
degree
of $\Cal K(M)$. Consider an $f\in\Cal K(M)$. Then $f$ is algebraically
dependent on $f_1,f_2,\hdots, f_d$; i.e. it satisfies an equation
\form{f^\lambda+g_1f^{\lambda-1}+\cdots+g_\lambda=0\, ,
}
\noindent
where $g_1,g_2,\hdots,g_\lambda\in\Bbb C(f_1,f_2,\hdots, f_d)$.
The minimal $\lambda$ for which such an equation holds is called the
{\it degree} of $f$ over $\Bbb C(f_1,f_2,\hdots, f_d)$. The content
of Proposition \teordb\; is that this degree is bounded from above by
$\kappa=\kappa(f_1,f_2,\hdots, f_d)$. Now choose an element
$\Theta\in\Cal K(M)$ so that its degree $\alpha$ is maximal. For any
$f\in\Cal K(M)$ consider the algebraic extension field
$\Bbb C(f_1,f_2,\hdots, f_d,\Theta,f)$. By the primitive element
theorem this extension is simple; i.e. there exists an element
$h\in\Bbb C(f_1,f_2,\hdots, f_d,\Theta,f)$ such that
$\Bbb C(f_1,f_2,\hdots, f_d,\Theta,f)= \Bbb C(f_1,f_2,\hdots, f_d,h)$.
Then 
\form{\matrix\format\r&\l\\
 \alpha&\qquad \geq \qquad
[\Bbb C(f_1,f_2,\hdots, f_d,h):\Bbb C(f_1,f_2,\hdots, f_d)]\\
\\
&=[\Bbb C(f_1,f_2,\hdots, f_d,\Theta,f):\Bbb C(f_1,f_2,\hdots, f_d,\Theta)]\\
&\qquad\qquad\qquad \times \,
[\Bbb C(f_1,f_2,\hdots, f_d,\Theta):\Bbb C(f_1,f_2,\hdots, f_d)]\\
\\
&\qquad \qquad \qquad\qquad \geq\alpha\, .\\
\endmatrix}
Hence the first factor on the right must be one; therefore
$f\in\Bbb C(f_1,f_2,\hdots, f_d,\Theta)$. The conclusion is that
\form{\Cal K(M)= \Bbb C(f_1,f_2,\hdots, f_d,\Theta)=
\Bbb C(f_1,f_2,\hdots, f_d)[\Theta]\, ,
}\edef\formudj{\number\q .\number\t}
\noindent
and {\sl any $f\in\Cal K(M)$ can be written as a polynomial of degree
$<\alpha$ having coefficients that are rational functions of
$f_1,f_2,\hdots, f_d$.}\par
From the above remark we derive the
\prop{There is an open neighborhood $U$ of $M$ in $X$ such that the
restriction map
\form{ \Cal K(U) @>>>\Cal K(M)}
\noindent
is an isomorphism.}
Here $\Cal K(U)$ denotes the field of meromorphic functions on $U$.
\smallskip
Let $M$ be a connected smooth abstract $CR$ manifold of type $(n,k)$.
Consider a complex $CR$ line bundle $F@>\pi>>M$ over $M$. Introduce
the graded ring
\form{\Cal A(M,F)=\bigcup_{\ell=0}^\infty{\CR(M,F^\ell)}\, ,
}
\noindent
where $\CR(M,F^\ell)$ are the smooth global $CR$ sections of the
$\ell$-th tensor power of $F$. Note that if $\sigma_1\in\CR(M,F^{\ell_1})$
and $\sigma_2\in\CR(M,F^{\ell_2})$, then
$\sigma_1\sigma_2\in\CR(M,F^{\ell_1+\ell_2})$.
\par
Assume that we are in a situation where smooth sections of $F$ have the weak
unique continuation property; e.g. we could take $M$ to be
essentially pseudoconcave (see [HN8]). Then
$\Cal A(M,F)$ is an integral domain because $M$ is connected. Let
\form{\Cal Q(M,F)=\left\{\dsize\frac{\sigma_1}{\sigma_2}\, \left|
\, \sigma_1,\sigma_2\in\CR(M,F^\ell)\;\text{for some $\ell$, and}\;
\sigma_0\not\equiv 0\right.\right\}}
\noindent
denote the field of quotients. \par
Then $\Cal Q(M,F)\subset\hat{\Cal K}(M)$,
and $\CR(M)=\Cal A(M,\text{trivial bundle})$.
\prop{Assume that $M$ is compact and has property $E$. 
\roster
\item 
If $F$ is locally
$CR$ trivializable, then $\Cal Q(M,F)$ is an algebraically closed subfield
of $\Cal K(M)$.
\item 
There exists a choice of a locally $CR$ trivializable $F$ such that
$\Cal Q(M,F)=\Cal K(M)$.
\par
Assume that $M$ is compact and essentially pseudoconcave. Then
\item $\Cal Q(M,F)$ is algebraically closed in $\hat{\Cal K}(M)$.
\par
In case $M$ is compact and satisfies both hypothesis, then
\item ${\Cal K}(M)$ is algebraically closed in $\hat{\Cal K}(M)$.
\endroster}
\demo{Proof}
To prove (1) [or (3)] we take an $h\in\Cal K(M)$ [or $h\in\hat{\Cal K}(M)$]
which is algebraic over $\Cal Q(M,F)$; i.e. $h$ satisfies an equation
of minimal degree
\form{h^\mu+k_1h^{\mu-1}+\cdots +k_\mu\, 
\equiv\, 0\, ,}\edef\fordn{\number\q.\number\t}\noindent
where $k_i\in\Cal Q(M,F)$. Let $k_i=\dsize\frac{s_i}{t_i}$ with
$s_i,t_i\in\CR(M,F^{\ell_i})$.
Multiplying by $\sigma_0=\prod_{i=1}^\mu{t_i}$ we obtain
\form{\sigma_0 h^\mu+\sigma_i h^{\mu-1}+\cdots +\sigma_\mu\,
\equiv \, 0\, ,}\noindent
where $\sigma_0\not\equiv 0$ and $\sigma_i\in\CR(M,F^\ell)$ for
$\ell=\sum_{i=1}^\mu{\ell_i}$. Multiplying by $\sigma_0^{\mu-1}$ we have:
\form{P(\sigma_0 h)=
\left(\sigma_0 h\right)^\mu+\sigma_1\left(\sigma_0 h\right)^{\mu-1}+
\cdots +\sigma_0^{\mu-2}\sigma_{\mu-2}\left(\sigma_0 h\right)
+\sigma_0^{\mu-1}\sigma_\mu
\equiv 0\, .}\edef\fordp{\number\q.\number\t}\noindent
Note that $\left(\sigma_0 h\right)$ is bounded, and hence is a smooth
section of $F^\ell$ over $M$. In a local smooth trivialization
of the bundle $F$, the tangential Cauchy-Riemann operator on sections
$s$ of $F$ has the form (see [HN8])
\form{\bar\partial^F_Ms=\bar\partial_M s+As\, 
,}\noindent
where $A$ is a smooth $\bar\partial_M$-closed $(0,1)$ form on $M$.
For the $\ell$-th tensor power $F^\ell$ of $F$, in the same trivialization,
we have
\form{\bar\partial^{F^\ell}_Ms=\bar\partial_M s+\ell As\, 
.}\noindent
We apply $\bar\partial^{F^\ell\mu}_M$ to both sides of (\fordp) and
obtain
\form{0=\bar\partial^{F^{\ell\mu}}_MP(\sigma_0h)=P'(\sigma_0h)\left[
\bar\partial_M(\sigma_0h)+\ell A\cdot(\sigma_0h)\right]=
P'(\sigma_0h)\bar\partial^{F^\ell}_M(\sigma_0h)\, 
.}\noindent
Because $\mu$ is minimal 
in (\fordn), we obtain that $P'(\sigma_0h)\not\equiv 0$,
and therefore $\tau=\sigma_0 h\in\CR(M,F^\ell)$. Since
$\sigma_0\in\CR(M,F^\ell)$, we obtain that $h=\dsize\frac{\tau}{\sigma_0}
\in\Cal Q(M,F)$. This completes the proof of (1) and (3).
\smallskip
To prove (2) it suffices to observe that it is possible to choose a
locally trivializable smooth complex $CR$ line bundle $F$ over $M$ such
that the $f_1,f_2,\hdots, f_d,\Theta$ appearing in (\formudj)
belong to $\Cal Q(M,F)$. Then by (\formudj), $\Cal Q(M,F)=\Cal K(M)$.
\smallskip
Finally (4) is a consequence of (2) and (3). This completes the proof
of the proposition.
\enddemo
Let $M$ be a connected smooth compact $CR$ manifold of type $(n,k)$,
having property $E$. Then
\thm{Let $F@>\pi>> M$ be a locally $\CR$ trivializable smooth
complex $CR$ line bundle over $M$. Then
\form{\roman{dim}_{\Bbb C}{\CR(M,F)}<\infty\, 
.}\noindent}
\demo{Proof}
For each point $a\in M$, we fix an open neighborhood $\omega_a$ of $a$ in $M$,
and an open coordinate neighborhood $\Omega_a$ of $a$ in $X$, such that
$\omega_a=\Omega_a\cap M$ and the restriction map
$\Cal O(\Omega_a)@>>>\CR(\omega_a)$ is an isomorphism. These $\omega_a$
are to be chosen to give a local $CR$ trivialization of $F$. Denote the
$\CR$ transition functions by $g_{ab}$. We may assume that the
$g_{ab}$ are bounded on $\omega_a\cap\omega_b$. Introduce the polydiscs
$K_a$ and $L_a$ as in (\formbe), and choose $a_1,a_2,\hdots,a_m$ on $M$
so that the $L_{a_1},L_{a_2},\hdots,L_{a_m}$ provide an open covering
of $M$. Choose an integer $\mu$ such that
\form{\left| g_{ab}\right|\, <\, e^\mu\quad\text{on}\quad \omega_a\cap
\omega_b}\noindent
for $a,b=a_1,a_2,\hdots,a_m$. Consider a section $s\in\CR(M,F)$
which vanishes at each point $a_i$ of order $\geq \mu+1$. The section
$s$ is represented by a smooth $CR$ function $s_i$ on each $\omega_{a_i}$.
Let
\form{\Upsilon=\max_{1\leq i\leq m}\max_{K_{a_i}}{|\tilde s_i|}\, ,
}\noindent
as before.
This maximum is attained at some point $z^*$ belonging to some
$K_{a^*}$, where $a^*$ is one of the $a_1,a_2,\hdots,a_m$.
Then by (iv) in Corollary 1.2, there is some $z^{**}\in\omega_{a^*}$ 
such that
\form{\tilde s_{a^*}(z^*)=s_{a^*}(z^{**})\, 
.}\noindent
But $z^{**}$ belongs to some $L_{a^{**}}$, where $a^{**}$ is one of the
$a_1,a_2,\hdots,a_m$. By the Schwarz lemma,
\form{\left|s_{a^{**}}(z^{**})\right|\,\leq\, \Upsilon\,e^{-(\mu+1)}\, 
.}\noindent
On the other hand we have
\form{\Upsilon=\left|s_{a^{*}}(z^{**})\right|=
\left| g_{a^*a^{**}}(z^{**})\right|\cdot 
\left|s_{a^{**}}(z^{**})\right|\, \leq \, M\, e^{-1}\, 
.}\noindent
This shows that the Taylor expansion of the representatives $s_i$, at $a_i$,
up to order $\mu$, completely determine a section $s$. Hence
$\roman{dim}_{\Bbb C}{\CR(M,F)}\leq m\binom{n+k+\mu+1}{\mu+1}<\infty$.
\enddemo
\se{The Chow theorem for $CR$ manifolds}
Again $M$ is a connected smooth compact $CR$ manifold of type $(n,k)$,
having property $E$. We have:
\thm{Let $\tau:M@>>>\Bbb C\Bbb P^N$ be a smooth $CR$ map. Suppose that
$\tau$ has maximal rank $2n+k$ at one point of $M$. Then
$\tau(M)$ is contained in an irreducible algebraic subvariety
of complex dimension $n+k$, and the transcendence degree of
$\Cal K(M)$ over $\Bbb C$ is $n+k$.}
\demo{Proof}
Let $Y$ be the smallest algebraic subvariety of $\Bbb C\Bbb P^N$ 
containing $\tau(M)$. Certainly $Y$ exists and is irreducible; it is
defined by the homogeneous prime ideal
\form{\Cal P_Y\, = \,
\{p\in\Bbb C_0[z_0,z_1,\hdots,z_N]\, | \, p\circ\tau=0\, \}\,
,}\edef\forea{\number\q.\number\t}\noindent
where $\Bbb C_0[z_0,z_1,\hdots,z_N]$ is the graded ring of
homogeneous polynomials on $\Bbb C\Bbb P^N$.\par
Let $\Cal R(Y)$ be the field of rational functions on $Y$. Any element
$f\in\Cal R(Y)$ is represented as the quotient of two homogeneous
polynomials $f=p/q$, with $q\notin\Cal P_Y$. If $f=p/q=p'/q'$, then
$pq'-p'q\in\Cal P_Y$. By (\forea) this shows that $f\circ\tau$ is
a well defined $CR$ meromorphic function in $\Cal K(M)$, and that the
homomorphism
\form{\tau^*:\Cal R(Y) @>>>\Cal K(M)
}\edef\foreb{\number\q.\number\t}\noindent
is injective. By the assumption that $\tau$ has maximum rank at one point,
and by Theorem 2.1, we obtain
\form{\matrix\format\r&\l\\
n+k\leq\roman{dim}_{\Bbb C}{Y}&=\text{transcendence degree of
$\Cal R(Y)$}\\
&\leq\text{transcendence degree of
$\Cal K(M)$}\leq n+k\, .
\endmatrix}\noindent
This completes the proof of the theorem.
\enddemo
As a corollary we obtain for compact
$CR$ manifolds an analogue of Chow's theorem:
\thm{Let $M$ be a connected smooth compact $CR$ embedded submanifold
of $\Bbb C\Bbb P^N$, of type $(n,k)$ and having property $E$.
Then $M$ is a generically embedded $CR$ submanifold of an irreducible
algebraic subvariety $Y$ of $\Bbb C\Bbb P^N$; moreover $M$
is contained in the set
$\roman{reg}\,{Y}$ of regular points of $Y$, 
and the map (\foreb) is an isomorphism.}\par
Setting $k=0$ above, we recover the theorem of Chow [C].
\demo{Proof}
We verify that $M$ avoids any singularity of $Y$. Consider a point
$x_0\in M$. We can assume that near $x_0$, the inhomogeneous coordinates
are $z_1,z_2,\hdots,z_N$, centered at $x_0$, and that $z_1,z_2,\hdots,
z_{n+k}$ are coordinates for the smallest affine complex linear subspace
containing $T_{x_0}M$. Let $\pi:\Bbb C^N@>>>\Bbb C^{n+k}$ be the
holomorphic projection onto this affine subspace. Then locally near
$x_0$, $M'=\pi(M)$ is a smooth generic $CR$ submanifold of $\Bbb C^{n+k}$,
having property $E$. But $\phi=\left(\pi\left|_{\dsize M}\right.\right)^{-1}$
is a $CR$ map on $M'$, so by property $E$ it has a local holomorphic
extension $\tilde\phi$ to a neighborhood $\omega$ of $\pi(x_0)$ in
$\Bbb C^{n+k}$. Then
$z_j-\tilde\phi_j(z_1,z_2,\hdots,z_{n+k})$, for $j=n+k+1,\hdots,N$, give
$N-(n+k)$ holomorphic functions in $\pi^{-1}(\omega)$ which vanish locally
on $M$ near $x_0$, and have independent differentials. In particular they
define germs in $\Cal P_Y\otimes\Cal O_{x_0}$. Using the fact that
$\Cal O_{x_0}$ is faithfully flat over the ring $\Bbb C[z_1,z_2,\hdots,z_N]$,
we obtain that $Y$ is locally a smooth complete intersection at $x_0$.
Thus $M\subset\roman{reg}\, Y$.
\par
Since (\foreb) is injective, $\Cal K(M)$ has a transcendence basis
$f_1,f_2,\hdots,f_{n+k}$ with each $f_i=\tilde f_i\left|_{\dsize M}\right.$,
where $f_i\in\Cal R(Y)$. Thus any $f\in\Cal K(M)$ is the solution of
an irreducible algebraic equation of the form
\form{f^\lambda+g_1 f^{\lambda-1}+\cdots + g_{\lambda}=0\, 
,}\edef\fored{\number\q.\number\t}\noindent
with coefficients $g_1,g_2,\hdots,g_\lambda\in\Cal R(Y)$. Take a point
$x_0\in M$ where $g_1,g_2,\hdots,g_\lambda$ are smooth $CR$ functions,
and
\form{\xi^\lambda+g_1 \xi^{\lambda-1}+\cdots + g_{\lambda}=0\, 
,}\edef\foree{\number\q.\number\t}\noindent
has $\lambda$ distinct complex roots. Since $\Cal R(Y)$ is algebraically
closed, there is some root $\tilde f\in\Cal R(Y)$ of (\fored) such
that $\tilde f(x_0)=f(x_0)$. As the roots of (\foree) depend continuously
on the coefficients, we have that $\tilde f(x)=f(x)$ for $x$ in
a neighborhood of $x_0$ on $M$. By unique continuation,
$f=\tilde f\left|_{\dsize M}\right.$. This completes the proof of the
theorem.\enddemo
\se{Projective embedding}
Let $M$ be a connected smooth compact $CR$ manifold of type $(n,k)$, having
property $E$. Then
\thm{The following are equivalent:
\roster
\item
$M$ has a smooth $CR$ embedding as a $CR$ submanifold of some 
$\Bbb C\Bbb P^N$.
\item There exists over $M$ a smooth complex $CR$ line bundle $F$ such that
the graded ring $\Cal A(M,F)=\bigcup_{\ell=0}^\infty{\CR(M,F^\ell)}$
separates points and gives ``local coordinates'' at each point of $M$.
\endroster}
\demo{Proof}
First we explain the meaning of (2). To say that $\Cal A(M,F)$
separates points means that, given $x\neq y$ on $M$, there exists
an integer $\ell=\ell(x,y)>0$, and smooth sections $s_0,s_1\in
\CR(M,F^\ell)$ such that
\form{\roman{det}\left[
\matrix
s_0(x)&s_0(y)\\
s_1(x)&s_1(y)
\endmatrix\right]\,\neq\, 0\, .}\edef\frmfa{\number\q.\number\t}
To say that $\Cal A(M,F)$
gives ``local coordinates'' means that, given $x$ in $M$, there 
exists an integer $\ell=\ell(x)>0$ 
and smooth sections $s_0,s_1,\hdots,s_{n+k}\in\CR(M,F^\ell)$
such that
\form{\dsize\sum_{i=0}^{n+k}{(-1)^is_i\,ds_0\wedge\cdots\wedge
\widehat{ds_i}\wedge\cdots\wedge ds_{n+k}}\neq 0\quad\text{at}\quad x\, 
.}\edef\frmfb{\number\q.\number\t}
In other words, the $CR$ meromorphic function 
$s_1/s_0$ takes {\it different values} at $x$ and $y$; and
(assuming $s_0(x)\neq 0$) the $CR$ meromorphic functions
$s_1/s_0$, $s_2/s_0$, $\hdots$, $s_{n+k}/s_0$ provide
analytically independent local $CR$ functions at $x$. Note that
a line bundle $F$ satisfying (2) is locally $CR$ trivializable; hence
also $F^\ell$.\smallskip
First we show that (2) implies (1): By row reduction we may assume that
the matrix (\frmfa) is in diagonal form; i.e. $s_0(y)=s_1(x)=0$ and
$s_0(x)\cdot s_1(y)\neq 0$. Hence, for any $m>0$,
\form{\roman{det}\left[
\matrix
s_0^m(x)&s_0^m(y)\\
s_1^m(x)&s_1^m(y)
\endmatrix\right]\,\neq\, 0\, .}
This shows that if two points $x$ and $y$ are separated by
$s_0,s_1\in\CR(M,F^\ell)$, they are also separated by
$s_0^m,s_1^m\in\CR(M,F^{\ell m})$.
Next we consider sections 
$s_0,s_1,\hdots,s_{n+k}\in\CR(M,F^\ell)$ which satisfy
(\frmfb) at a given point $x\in M$. We can assume $s_0(x)\neq 0$.
By changing $s_i$ to $s_i-(s_i(x)/s_0(x))s_0$, we can also assume that
$s_i(x)=0$ for $i=1,2,\hdots,n+k$. Then for any $m>0$, we have that
$s_0^m,\, s_0^{m-1}s_1,\,\hdots,\,s_0^{m-1}s_{n+k}\in\CR(M,F^{\ell m})$
also satisfy(\frmfb). For each $\ell>0$ we define the following
subset of $M\times M$:
\form{\matrix\format\r&\l\\
\Cal U_\ell&=\{(x,y)\, | \, \exists (s_0,s_1)\in\CR(M,F^\ell)\;
\text{such that (\frmfa) holds}\}\\
&\qquad\cup\,
\{(x,x)\, | \, \exists s_0,s_1,\hdots,s_{n+k}\in\CR(M,F^\ell)\; 
\text{such that (\frmfb) holds}\}\,\}
\endmatrix}
Each $\Cal U_\ell$ is open in $M\times M$, and 
$\Cal U_\ell\subset \Cal U_{\ell\cdot m}$ for $m>0$. By our hypothesis
(2) the $\{\Cal U_\ell\}$ give an open covering of
$M\times M$. By compactness $M\times M$ is covered by
$\Cal U_{\ell_1},\, \Cal U_{\ell_2}\,\hdots,\, \Cal U_{\ell_r}$.
Therefore $M\times M=\Cal U_{\ell_0}$ where $\ell_0=\ell_1\cdot\ell_2\cdots
\ell_r$.
Now let $s_0,s_1,\hdots,s_N$ be a basis for smooth sections of the
complex vector space $\CR(M,F^{\ell_0})$, which is finite dimensional
by Theorem 4.5. Consider the map
$\tau:M@>>>\Bbb C\Bbb P^N$
defined by
$x@>>>(s_0(x):s_1(x):\hdots:s_N(x))$.
It is a $CR$ map, which is one-to-one, and an immersion at each point;
thereby giving a $CR$ embedding of $M$ into $\Bbb C\Bbb P^N$.
\smallskip
Next we show that (1) implies (2): It suffices to take $F$ equal to the
pull-back to $M$ of the hyperplane section line bundle on $\Bbb C\Bbb P^N$.
With this $F$ we may take $\ell=1$ in (\frmfa) and (\frmfb).
\par
This completes the proof.
\enddemo
\vfill
\pagebreak


\Refs\widestnumber\key{aaaaa}

\hyphenation{con-ca-ve}
\ref\key A
\by A.Andreotti
\paper Th\'eor\`emes de dependence alg\'ebrique sur les espaces
complexes pseudo-concaves
\jour Bull. Soc. Math. France
\vol 91
\yr 1963
\pages 1-38
\endref

\ref\key AG
\by A.Andreotti, H.Grauert
\paper Algbebraische K\"orper von automorphen Funktionen
\jour Nachr. Ak. Wiss. G\"ottingen
\yr 1961
\pages 39-48
\endref



\ref\key BHN
\by J.Brinkschulte, C.D.Hill, M.Nacinovich
\paper Remarks on weakly pseudoconvex boundaries
\pages 1-9
\yr 2001
\jour Preprint
\moreref Indagationes Mathematicae \toappear
\endref

\ref\key BP
\by A.Boggess, J.Polking
\paper Holomorphic extensions of $CR$ functions
\jour Duke Math. J.
\vol 49
\yr 1982
\pages 757-784
\endref

\ref\key C
\by W.L.Chow
\paper On complex compact analytic varieties
\jour Amer. J. Math
\vol 71
\yr 1949
\pages 893-914
\endref


\ref\key DCN
\by L.De Carli,
M.Nacinovich
\paper Unique continuation in abstract $CR$ manifolds
\jour
Ann. Scuola Norm. Sup. Pisa Cl. Sci.
\vol 27
\yr 1999
\pages
27-46
\endref

\ref \key HN1
\by C.D.Hill, M.Nacinovich
\paper A necessary condition for global Stein immersion of compact
$CR$ manifolds
\jour Riv. Mat. Univ. Parma
\vol 5
\yr
1992
\pages 175-182
\endref

\ref \key HN2
\bysame
\paper The topology of Stein CR manifolds and the Lefschetz theorem
\jour Ann. Inst.
Fourier, Grenoble
\vol 43
\yr 1993
\pages 459--468
\endref


\ref \key HN3
\bysame 
\paper Pseudoconcave $CR$ manifolds
\inbook Complex Analysis and
Geometry, (eds Ancona, Ballico, Silva)
\publ Marcel Dekker, Inc
\yr
1996
\publaddr New York
\pages 275--297
\endref


\ref \key HN4
\bysame
\paper Aneurysms of pseudoconcave $CR$ manifolds
\jour Math.
Z.
\vol 220
\pages 347--367
\yr 1995
\endref

\ref \key HN5
\bysame
\paper
Duality and distribution cohomology of $CR$ manifolds
\jour Ann. Scuola
Norm. Sup. Pisa
\vol 22
\yr 1995
\pages 315--339
\endref


\ref \key HN6
\bysame
\paper On the Cauchy problem in complex analysis
\jour Annali
di matematica pura e applicata
\vol CLXXI (IV)
\yr 1996
\pages
159-179
\endref


\ref\key HN7
\bysame
\paper Conormal
suspensions of differential complexes
\jour J. Geom. Anal. 
\vol 10
\yr
2000
\pages 481-523
\endref

\hyphenation{qua-der-ni}
\ref\key HN8
\bysame
\paper A weak
pseudoconcavity condition for abstract almost $CR$ manifolds
\jour Invent.
math.
\yr 2000
\vol 142
\pages 251-283
\endref

\ref\key HN9
\bysame
\paper Weak pseudoconcavity and the maximum modulus principle
\jour Quaderni sez. Geometria Dip. Matematica Pisa
\toappear
\; in Ann. Mat. Pura e Appl.
\pages 1-10
\yr 2001
\endref

\ref\key HN10
\bysame
\paper Pseudoconcavity at infinity
\jour Quaderni sez. Geometria Dip.
Matematica Pisa
\yr  2000
\vol 1.229.1228
\pages 1-27
\endref


\ref\key HN11
\bysame
\paper Two lemmas on double complexes and their applications to
$CR$ cohomology
\jour Quaderni sez. Geometria Dip.
Matematica Pisa
\yr 2001
\pages 1-10
\vol 1.262.1331
\endref

\ref \key L
\by E.E.Levi
\paper Studii sui punti singolari essenziali delle funzioni analitiche
di due o pi\`u variabili complesse
\jour Ann. Mat. Pura Appl.
\vol XVII (s.III)
\yr 1909
\moreref Opere, Cremonese, Roma, 1958, 187-213
\endref

\ref \key NV
\by
M.Nacinovich, G.Valli
\paper Tangential Cauchy-Riemann complexes on
distributions
\jour Ann. Mat. Pura Appl.
\vol 146
\yr 1987
\pages
123-160
\endref

\ref \key Se
\by J.P.Serre
\paper Fonctions automorphes, quelques majorations dans le cas o\`u
$X/G$ est compact.
\inbook S\'eminaire H. Cartan 1953-54
\yr 1957
\publ Benjamin
\publaddr New York
\endref

\ref\key Si
\by C.L.Siegel
\paper Meromorphe Funktionen auf kompakten analytischen
Manningfaltigkeiten
\jour Nachr. Ak. Wiss G\"ottingen
\yr 1955
\pages 71-77
\endref
\endRefs
\enddocument